\documentclass[12pt]{gtart}

\usepackage{amsmath, amssymb, latexsym} 
\usepackage{euscript}

\usepackage{mathptmx}

\def\bkR{{\rm I\kern-.17em R}}
\def\R{\bkR}
\def\bkH{{\rm I\kern-.17em H}}
\def\H{\bkH}

\newcommand{\Stab}{\operatorname{Stab}}
\newcommand{\Isom}{\operatorname{Isom}}
\def\Iso{\Isom^+(\H^n)}
\def\cprime{\/{\mathsurround=0pt$'$}}

\theoremstyle{plain}

\newtheorem{thm}{Theorem}
\newtheorem*{thm*}{Theorem}
\newtheorem{lem}[thm]{Lemma}
\newtheorem*{lem*}{Lemma}
\newtheorem{cor}[thm]{Corollary}
\newtheorem*{cor*}{Corollary}
\newtheorem{pro}[thm]{Proposition}
\newtheorem*{pro*}{Proposition}

\newtheorem*{rem*}{Remark}

\newtheorem*{defn*}{Definition}

\newtheorem*{ackn}{Acknowledgements}

\begin{document}
\title{Geodesic ideal triangulations exist virtually}
\author{Feng Luo, Saul Schleimer and Stephan Tillmann}

\begin{abstract}
It is shown that every non-compact hyperbolic manifold of finite volume has a finite cover admitting a geodesic ideal triangulation. Also, every hyperbolic manifold of finite volume with non-empty, totally geodesic boundary has a finite regular cover which has a geodesic partially truncated triangulation. The proofs use an extension of a result due to Long and Niblo concerning the separability of peripheral subgroups.
\end{abstract}


\maketitle
\footnotetext{This work is in the public domain.}


Epstein and Penner~\cite{EpPe1988} used a convex hull construction in Lorentzian space to show that every non-compact hyperbolic manifold of finite volume has a canonical subdivision into convex geodesic polyhedra all of whose vertices lie on the sphere at infinity of hyperbolic space. In general, one cannot expect to further subdivide these polyhedra into ideal geodesic simplices such that the result is an ideal triangulation. That this is possible after lifting the cell decomposition to an appropriate finite cover is the first main result of this paper. A cell decomposition of a hyperbolic $n$--manifold into ideal geodesic $n$--simplices all of which are embedded will be referred to as an \emph{embedded geodesic ideal triangulation.}

\begin{thm}\label{thm:ideal triangulation}
Any non-compact hyperbolic manifold of finite volume has a finite regular cover which admits an embedded geodesic ideal triangulation.
\end{thm}

The study of geodesic ideal triangulations of hyperbolic 3--manifolds goes back to Thurston \cite{Th1978}. They are known to have nice properties through, for instance, work by Neumann and Zagier \cite{NeZa1985} and Choi \cite{Ch2004}. Petronio and Porti \cite{PePo} discuss the question of whether every non-compact hyperbolic 3--manifold of finite volume has a geodesic ideal triangulation --- this question still remains unanswered.

Kojima~\cite{Ko1990} extended the construction by Epstein and Penner to obtain a canonical decomposition into partially truncated polyhedra of any hyperbolic manifold with totally geodesic boundary components. A cell decomposition of a hyperbolic $n$--manifold with totally geodesic boundary into geodesic partially truncated $n$--simplices all of which are embedded will be referred to as an \emph{embedded geodesic partially truncated triangulation.}

\begin{thm}\label{thm:partially truncated triangulation}
Any finite-volume hyperbolic manifold with non-empty, totally geodesic boundary has a finite regular cover which admits an embedded geodesic partially truncated triangulation.
\end{thm}

An ideal polyhedron will be viewed as a special instance of a partially truncated one, which allows a unified proof of Theorems \ref{thm:ideal triangulation} and \ref{thm:partially truncated triangulation}. They are proved by showing that any cell decomposition lifts to some finite regular cover where it can be subdivided consistently. In particular, one has the following application. Kojima \cite{Ko1992} shows that any 3--dimensional finite-volume hyperbolic manifold with non-empty, totally geodesic boundary has a decomposition into geodesic partially truncated polyhedra each of which has at most one ideal vertex. Frigerio \cite{Fr2005} conjectures that such a decomposition exists where each polyhedron is a tetrahedron; a virtually affirmative answer is an immediate consequence of the proof of Theorem \ref{thm:partially truncated triangulation}:

\begin{cor}
Any 3--dimensional finite-volume hyperbolic manifold with non-empty, totally geodesic boundary has a finite regular cover which admits a decomposition into partially truncated geodesic tetrahedra each of which has at most one ideal vertex.
\end{cor}

The key result used in the proof of Theorem \ref{thm:partially truncated triangulation} is the following theorem which follows easily from work by Long and Niblo \cite{LoNi1991}. A subgroup $H$ of a group $G$ is \emph{separable} in $G$ if given any element $\gamma \in G \setminus H,$ there is a finite index subgroup $K \le G$ which contains $H$ but $g \notin K.$ If $M$ is a hyperbolic manifold of finite volume with (possibly empty) totally geodesic boundary, then a subgroup of $\pi_1(M)$ is termed \emph{peripheral} if it is either conjugate to the fundamental group of a totally geodesic boundary component or to the fundamental group of a cusp or $\partial$--cusp.

\begin{thm}[Long--Niblo]\label{thm:Long-Niblo}
Let $M$ be a hyperbolic manifold of finite volume with (possibly empty) totally geodesic boundary. Then every peripheral subgroup of $M$ is separable in $\pi_1(M).$
\end{thm}


\begin{ackn}
\emph{The authors thank Chris Leininger for an enlightening conversation. They also thank Daryl Cooper for bringing the work by Long and Niblo to their attention. Research of the first author is supported in part by the NSF. The second author is partly supported by the NSF (DMS-0508971). The third author is supported under the Australian Research Council's Discovery funding scheme (project number DP0664276).}
\end{ackn}


\section{Subgroup separability}
\label{sec:Subgroup separability}

Let $M$ be a finite-volume hyperbolic $n$--manifold with non-empty totally geodesic boundary. Following Kojima \cite{Ko1990}, the periphery of $M$ is made up of three parts:
first, $\partial M$ consisting of totally geodesic closed or non-compact hyperbolic $(n-1)$--manifolds;
second, \emph{(internal) cusps} modelled on closed Euclidean $(n-1)$--manifolds;
and third, \emph{$\partial$--cusps} modelled on compact Euclidean $(n-1)$--manifolds with geodesic boundary. The boundary of  $\partial$--cusps is contained on non-compact geodesic boundary components.

For the remainder of this paper, $M$ denotes a finite-volume hyperbolic manifold which is either non-compact or has non-empty totally geodesic boundary. Without loss of generality, it may be assumed that $M$ is orientable. Note that either $\widetilde{M} = \H^n$ or it can be viewed as the complement of an infinite set of hyperplanes in $\H^n;$ in either case there is an identification $\pi_1(M) = \Gamma \le \Iso.$ 

\begin{pro}[Long--Niblo]
Let $X$ be a totally geodesic component of $\partial M.$ Choose a basepoint $x \in X.$ Then $\pi_1(X,x)$ is a separable subgroup of $\pi_1(M,x).$
\end{pro}

\begin{proof}
Let $D$ denote the manifold obtained by doubling $M$ along $X.$ Then $D$ is hyperbolic with (possibly empty) totally geodesic boundary, and hence $\pi_1(D)  \le \Iso$ is residually finite due to a result by Mal\cprime cev \cite{Ma1940}. The proof in \S 2 of Long and Niblo \cite{LoNi1991} now applies to this set-up.
\end{proof}

\begin{pro}
Let $X$ be a horospherical cross section of a cusp or $\partial$--cusp of $M.$ Choose a basepoint $x \in X.$ Then $\pi_1(X,x)$ is a separable subgroup of $\pi_1(M,x).$
\end{pro}

\begin{proof}
This follows from the well-known result that a maximal abelian subgroup of a residually finite group $\Gamma$ is separable in $\Gamma$ (see Ratcliffe \cite{Ra}).
\end{proof}

\begin{proof}[Proof of Theorem \ref{thm:Long-Niblo} ]
First note that if $\varphi\co G_1 \to G_2$ is an isomorphism and $H \le G_1$ is separable in $G_1,$ then $\varphi(H) \le G_2$ is separable in $G_2.$ In particular, reference to base points can be omitted. Next, note that if $H\le \Gamma$ is separable, so is $\gamma^{-1}H\gamma$ for any $\gamma \in \Gamma.$ Thus, Theorem \ref{thm:Long-Niblo} follows from the above propositions for orientable manifolds. If $M$ is non-orientable, denote by $\Gamma_0$ a subgroup of index two of $\Gamma$ corresponding to the fundamental group of the orientable double cover. Then any subgroup $H \le \Gamma$ is separable in $\Gamma$ if and only if $H \cap \Gamma_0$ is separable in $\Gamma_0.$ Now if $H$ is a peripheral subgroup of $\Gamma,$ then $H \cap \Gamma_0$ is a peripheral subgroup of $\Gamma_0.$
\end{proof}


\section{Partially truncated polyhedra}
\label{sec:Partially truncated polyhedra}

Certain convex geodesic polyhedra in $\H^n$ are termed \emph{geodesic partially truncated polyhedra} and can be described intrinsically. However, reference to the projective ball model $B^n \subset \R^n$ will be made here, and $\H^n$ will be identified with $B^n$ throughout.

Let $\hat{P}$ be an $n$--dimensional convex Euclidean polyhedron in $\R^n$ such that (1) each vertex is either called \emph{ideal} or \emph{hyperideal,} (2) its ideal vertices are contained on $\partial B^n,$ (3) its hyperideal vertices are contained in $\R^n \setminus \overline{B}^n,$ and (4) each face of codimension two meets $\overline{B}^n.$ Then a convex geodesic polyhedron $P \subset B^n$ is obtained by truncating $\hat{P}$ along hyperplanes canonically associated to its hyperideal vertices as follows. If $v \in \R^n$ is a hyperideal vertex then the associated hyperplane $H(v)$ is the hyperplane parallel to the orthogonal complement of $v$ (with respect to the standard Euclidean inner product on $\R^n$) which meets $\partial B^n$ in the set of all points $x$ with the property that there is a tangent line to $\partial B^n$ passing through $x$ and $v.$ The polyhedron $P$ is termed a \emph{geodesic partially truncated polyhedron,} and $\hat{P}$ the \emph{Euclidean fellow} of $P.$ Combinatorially, $P$ is obtained from $\hat{P}$ by removing disjoint open stars of all the hyperideal vertices as well as all the ideal vertices.

If a codimension one face of $P$ is contained in a face of $\hat{P},$ then it is called \emph{lateral}; otherwise it is a \emph{truncation face}. Lateral faces and truncation faces meet at right angles. If $P$ has no truncation faces, then it is also termed a \emph{geodesic ideal polyhedron.} 

Any subdivision of $\hat{P}$ into $n$--simplices without introducing new vertices uniquely determines a subdivision of $P$ into \emph{geodesic partially truncated $n$--simplices.} The polyhedron $\hat{P}$ is termed the \emph{Euclidean fellow} of $P.$


\section{The pulling construction}
\label{sec:The pulling construction}

Let $(\mathcal{C}, \Phi)$ be a geodesic partially truncated cell decomposition of $M,$ that is, $\mathcal{C}$ is a disjoint union of geodesic partially truncated polyhedra, each element in $\Phi$ is an isometric face pairing, and $M = \mathcal{C}/\Phi.$ Then $(\mathcal{C}, \Phi)$ pulls back to a $\Gamma$--equivariant cell decomposition of $\widetilde{M} \subseteq B^n,$ and for each $P \in \mathcal{C}$ one may choose an isometric lift $\tilde{P}$ to $B^n$ and hence a Euclidean fellow $\hat{P}\subset \R^n.$ The hyperideal vertices of $\hat{P}$ correspond to totally geodesic boundary components of $M,$ the ideal vertices of $\hat{P}$ to internal cusps, and the intersection of codimension-two faces of $\hat{P}$ with $\partial{B}^n$ to $\partial$--cusps.

Let $\hat{\mathcal{C}}= \cup \{\hat{P}\}$ be the finite disjoint union of the Euclidean fellows, and view $P \subset \hat{P}.$ The cell decomposition of $M$ induces face pairings $\hat{\Phi}$ such that $M$ is obtained from the pseudo-manifold $\hat{M} = \hat{\mathcal{C}}/\hat{\Phi}$ by deleting the ideal vertices and open stars of the hyperideal vertices, and each element of $\hat{\Phi}$ restricts to an element of $\Phi.$

\begin{lem}[Pulling construction]\label{lem:Pulling construction}
Suppose that no polyhedron in $\hat{\mathcal{C}}$ has two distinct vertices identified in $\hat{M}.$ Then $M$ has an embedded geodesic partially truncated triangulation.
\end{lem}

\begin{proof}
It suffices to show that there is a subdivision of $(\hat{\mathcal{C}}, \hat{\Phi})$ such that (1) each polyhedron in $\hat{\mathcal{C}}$ is simplicially subdivided into straight Euclidean $n$--simplices without introducing new vertices, and (2) the elements of $\hat{\Phi}$ restrict to simplicial face pairings with respect to the subdivision.

Choose an ordering of the cusps and totally geodesic boundary components of $M.$ This determines a well-defined, unique ordering of the 0--skeleton of $\hat{M}$ and, by assumption, of the vertices of each polyhedron in $\hat{\mathcal{C}}.$ One thus obtains the following unique subdivision of each polyhedron.

Let $P \in \mathcal{C},$ and label its vertices $v_0, v_1, ..., v_k$ such that $v_i > v_j$ if $i<j.$ Subdivide $\hat{P}$ by coning to $v_0$ each element of its $i$--skeleton, $0 \le i \le n-1,$ which does not contain $v_0.$ The result is a collection of polyhedra, $\mathcal{P}_0,$ together with well-defined face pairings $\Phi_0$ such that the identification space $\mathcal{P}_0/\Phi_0$ is $\hat{P}.$ One now proceeds inductively. Given $\mathcal{P}_j$ and $\Phi_j,$ subdivide each polyhedron in $\mathcal{P}_j$ containing $v_{j+1}$ by coning to $v_{j+1}$ each element of its $i$--skeleton, $0 \le i \le n-1,$ which does not contain $v_{j+1}.$ This gives $\mathcal{P}_{j+1},$ together with well-defined face pairings $\Phi_{j+1}$ such that the identification space $\mathcal{P}_{j+1}/\Phi_{j+1}$ is $\hat{P}.$

It needs to be shown that the set $\mathcal{P}_k$ is a collection of $n$--simplices. Indeed, let $Q \in \mathcal{P}_k,$ and assume that $v_h$ is its smallest vertex. Then $Q$ is the cone to $v_h$ of an $(n-1)$--dimensional face $F^{n-1}$ not containing $v_h.$ The face $F^{n-1}$ is the cone to its smallest vertex of an 
$(n-2)$--dimensional face $F^{n-2}$ not containing that vertex, and it follows inductively that $Q$ has exactly $n+1$ vertices.

Let $\hat{P}, \hat{P}' \in \hat{\mathcal{C}}$ with top-dimensional faces $\hat{F}$, $\hat{F}'$ such that there is a face pairing $\varphi \in \hat{\Phi}$ with $\varphi(\hat{F}) = \hat{F}'.$ The respective subdivisions of $\hat{F}$ and $\hat{F}'$ into ideal $(n-1)$--simplices depend uniquely on the ordering of their vertices. Whence $\varphi$ is simplicial with respect to the subdivisions, and restricts to a simplicial face pairing for each $n$--simplex in the subdivision. Moreover, the resulting decomposition of $\hat{M}$ is simplicial since any $n$--simplex has no two vertices identified, and hence must be embedded in $\hat{M}.$
\end{proof}


\section{Proof of the main results}
\label{sec:Proof of the main results}

The strategy of the proof is to create a finite regular cover $N$ of $M$ with the property that Lemma~\ref{lem:Pulling construction} can be applied to the pull back of $\mathcal{C}.$ The notation of the previous sections will be used. Recall that for each $P \in \mathcal{C},$ there is the fixed Euclidean fellow $\hat{P} \subset \R^n.$ The action of $\Gamma$ on $\overline{B}^n$ extends to $\R^n\setminus\overline{B}^n$ via the action on the associated hyperplanes. In particular, if $v \in \hat{P}$ is a vertex, then the subgroup $\Stab_\Gamma(v) \le \Gamma$ is peripheral.

Let $\mathcal{D}(M)$ be the following set of pairs of points in $\R^n:$ $(v,w) \in \mathcal{D}(M)$ if and only if there is some $P \in \mathcal{C}$ such that $v$ and $w$ are distinct vertices of $\hat{P}.$ Note that $\mathcal{D}(M)$ is finite; its elements are termed \emph{diagonals} for $M.$ A diagonal $(v,w)$ is said to be \emph{returning} if there is $\gamma \in \Gamma$ such that $\gamma v = w.$ Note that the pulling construction can be applied to $\hat{\mathcal{C}}$ if no diagonal is returning. Hence assume that this is not the case.

If $p \co N \to M$ is a finite cover, then the cell decomposition $(\mathcal{C}, \Phi)$ pulls back to a cell decomposition of $N,$ and there is a corresponding set of diagonals for $N.$ If $P \in \mathcal{C}$ pulls back to $P_1,...,P_k,$ then (up to relabelling) one may choose $\hat{P}=\hat{P}_1.$ In particular, it may be assumed that $\mathcal{D}(M) \subset \mathcal{D}(N);$ any other element of $\mathcal{D}(N)$ is of the form $\gamma \cdot (v,w) = (\gamma v, \gamma w)$ for some $\gamma \in \Gamma$ and $(v,w) \in \mathcal{D}(M).$ This choice will be made throughout. If $(v,w)\in \mathcal{D}(M)$ is not a returning diagonal for $M,$ then it is also not a returning diagonal for $N.$

Assume that $(v,w) \in \mathcal{D}(M)$ is a returning diagonal.

\begin{lem}
There is a finite (possibly not regular) cover $p \co N_{(v,w)} \to M$ such that ${(v,w)} \in \mathcal{D}(N)$ is not a returning diagonal.
\end{lem}

\begin{proof}
Since ${(v,w)}$ is a returning diagonal for $M,$ there is $\gamma \in \Gamma$ such that $\gamma v = w.$ In particular, $\gamma \notin \Stab_\Gamma(v)$ because $v$ and $w$ are distinct. Since $\Stab_\Gamma(v)$ is a peripheral subgroup, Theorem \ref{thm:Long-Niblo} yields a finite index subgroup $K \le \Gamma$ which contains $\Stab_\Gamma(v)$ but $\gamma \notin K.$ Denote by $p \co N_{(v,w)} \to M$ the finite cover corresponding to the subgroup $K,$ i.e.\thinspace $N_{(v,w)} = \widetilde{M}/K.$

Assume that $(v,w) \in \mathcal{D}(N)$ is a returning diagonal. Then there is $\delta \in K$ with the property that $\delta v = w.$ Thus, $\gamma^{-1}\delta \in \Stab_\Gamma(v) \le K$ which implies $\gamma \in K.$ But this contradicts the choice of $K.$ Whence $(v,w)$ is not a returning diagonal for $N.$
\end{proof}

\begin{lem}
If $N \to M$ is a regular cover which factors through $N_{(v,w)},$ then no element of the orbit $\Gamma \cdot {(v,w)}$ can be a returning diagonal for $N.$
\end{lem}

\begin{proof}
If $N \to M$ is a cover which factors through $N_{(v,w)},$ then ${(v,w)}$ cannot be a returning diagonal for $N.$ If $N \to M$ is a regular cover, then no element of the orbit $\Gamma \cdot{(v,w)}$ is a returning diagonal for $N$ since the action of the group of deck transformations is transitive and $\pi_1(N)$ corresponds to a normal subgroup of $\pi_1(M).$
\end{proof}

For each diagonal $(v,w)$ choose a finite cover $N_{(v,w)} \to M$ with the property that ${(v,w)}$ is not a returning diagonal for $N_{(v,w)} .$ This gives a finite collection of covers, and one may pass to a common finite cover $N$ with the property that $N\to M$ is regular. In particular, no element of the orbit of any diagonal for $M$ can be a returning diagonal for $N.$ The cell decomposition $(\mathcal{C}, \Phi)$ of $M$ lifts to a polyhedral cell decomposition of $N,$ to which the pulling construction can thus be applied. This completes the proof of Theorems~\ref{thm:ideal triangulation} and \ref{thm:partially truncated triangulation}. \qed



\address{Department of Mathematics, Rutgers University, New Brunswick, NJ 08854, USA}
\address{Department of Mathematics, Rutgers University, New Brunswick, NJ 08854, USA}
\address{Department of Mathematics and Statistics, The University of Melbourne, VIC 3010, Australia}
\email{fluo@math.rutgers.edu, saulsch@math.rutgers.edu, tillmann@ms.unimelb.edu.au}
\Addresses


\end{document}